\documentclass[12pt,leqno]{article}
\usepackage[francais]{babel}
\usepackage[latin1]{inputenc}
\usepackage[T1]{fontenc}
\usepackage{amsmath,amsfonts,amssymb,amsthm,amscd}
\usepackage{bm}
\date{}

\newcommand\rg{\rightarrow}

\newcommand\U{\mathbb{U}}

\newcommand\C{\mathbb{C}}

\newcommand\N{\mathbb{N}}

\newcommand\R{\mathbb{R}}

\newcommand\Q{\mathbb{Q}}

\newcommand{\ld}{_{|\ell^2}}

\newcommand{\newsmile}{\smallsmile_{\hspace{-3mm}\displaystyle\smallfrown}}

\def\adots{\mathinner{\mkern1mu\raise1pt\vbox{\kern7pt\hbox{.}}
\mkern2mu\raise4pt\hbox{.}
\mkern2mu\raise7pt\hbox{.}\mkern1mu}}

\newtheorem{monlem}{\textsc{Lemme}}[section]

\begin{document}

\title{Sur le spectre des opérateurs rigides}

\author{Pierre Mazet et Eric Saias}

\maketitle


\begin{center}
Pour Mustapha Krazem,

à l'occasion de son soixantième anniversaire
\end{center}

\section{Introduction}

Soit $E$ un espace vectoriel normé sur $\C$. On appellera ici opérateur de $E$ toute application linéaire continue $u$ de $E$ dans $E$. Il est dit rigide (respectivement uniformément rigide) quand il existe une suite strictement croissante $(n_k)_{k\ge1}$ d'entiers telle que pour tout $x$ de $E$, on a $\lim\limits_{k\rg +\infty} u^{n_k}x=x$ (resp. telle que $\lim\limits_{k\rg +\infty} u^{n_k}=id_E$ dans l'espace vectoriel normé des opérateurs de $E$).

L'article \cite{CoMaPa} de Costakis, Manoussos et Parissis de 2014 est un travail de référence pour les opérateurs rigides et uniformément rigides. Ils passent en revue des familles classiques d'opérateurs (certaines familles d'opérateurs diagonaux, certaines familles d'opérateurs de composition, etc) et pour chacune d'entre elles, ils donnent un  critère de rigidité, voire d'uniforme rigidité. Ils se demandent en particulier à quelle condition nécessaire et suffisante un endomorphisme de dimension finie est rigide, et de même pour uniformément rigide.

Poser la question, c'est y répondre. Encore faut il la poser !
C'est tout le mérite de Costakis, Manoussos et Parissis de l'avoir fait. 

Cela les amène à effectuer une démarche en trois temps.
En notant $\U=\{z\in\C :|z|=1\}$ et $\sigma(u)$ le spectre de $u$, ils donnent donc d'abord la preuve (voir le theorem 4.1 et la remark 4.3 de \cite{CoMaPa}) du


\eject

\textsc{Théorème} CMP.

\textit{Soit $u$ un endomorphisme d'un $\C$ espace vectoriel de dimension finie. Les conditions suivantes sont équivalentes.}

\begin{tabular}{ll}
\kern-0.7cm$(i)$ &$u$\textit{ est rigide}\\

\kern-0.7cm$(ii)$ &$u$ \textit{est uniformément rigide}\\

\kern-0.7cm$(iii)$ &$u$ \textit{est diagonalisable et} $\sigma(u) \subset \U$.
\end{tabular}

\vspace{2mm}

Ils se demandent ensuite dans quelle mesure ce résultat s'étend à la dimension infinie. Dans cette direction, ils montrent (voir la proposition 2.19 de \cite{CoMaPa}) la

\vspace{2mm}

\textsc{Proposition} CMP.

\textit{Soit $u$ un opérateur uniformément rigide d'un espace de Banach complexe. Alors $\sigma(u) \subset \U$.}

\vspace{2mm}

(Nous donnons en Annexe 1 une preuve de ce résultat plus rapide que celle de \cite{CoMaPa}).

Enfin quid des opérateurs rigides ? La proposition CMP s'étend--elle à ceux--ci ? Dans leur question 2.20 ils posent deux questions précises.

\vspace{2mm}

\textbf{Question 1.}

Est--ce que tous les opérateurs rigides sont inversibles ?

\vspace{2mm}

\textbf{Question 2.}

Est--ce que pour tout opérateur rigide et inversible, son inverse est rigide ?

\vspace{2mm}

La réponse aux deux questions est négative. Nous introduisons au chapitre suivant la notion d'opérateur diagonal par blocs qui étend celle d'opérateur diagonal de $\ell^2(\N^*)$. On a alors le

\vspace{2mm}

\textsc{Théorème.}

\textit{Pour tout $r\in [0,1]$, il existe un opérateur diagonal par blocs et rigide dont le spectre est $\{\lambda\in \C : r \le |\lambda| \le 1\}$.}

\vspace{2mm}

{Pour $r=0$, un tel opérateur est non inversible. Cela répond donc négativement à la Question~$1$. Signalons qu'après avoir prouvé notre théorème, nous avons été informés par Sophie Grivaux que l'existence d'opérateurs rigides dont le spectre est le disque unité fermé peut également être établie à partir du travail \cite{GrRo}.}

{Pour $r\in ]0,1[$, l'inverse d'un tel opérateur n'est pas rigide. En effet son spectre est $\{\lambda\in \C : 1\le |\lambda| \le r^{-1}\}$. Or il est bien connu que le spectre d'un opérateur rigide est inclus dans le disque unité fermé (voir l'Annexe~$2$ pour une preuve rapide qui utilise le théorème de Banach--Steinhaus). Cela répond donc négativement à la Question~$2$.}

\vspace{2mm}

Nous avons donc utilisé que le rayon spectral d'un opérateur rigide est inférieur ou égal à un. En réalité, il est égal à un. En effet, pour tout opérateur rigide $u$ d'un Banach, 

\vspace{2mm}

\noindent(1.1) toute composante connexe du spectre de $u$ rencontre le cercle unité. 

\vspace{2mm}

Cela découle facilement du théorème de décomposition de Riesz (theorem D.2.1 de l'appendice D.2 de \cite{BaMa}) comme l'ont montré Costakis, Manoussos et Parissis à la proposition 2.11 de \cite{CoMaPa}, sous une hypothèse plus faible que la rigidité. En réalité, la même argumentation avait déjà été utilisée par Kitai en 1982 pour établir que les opérateurs $u$ d'un Banach qui satisfont la propriété

\vspace{2mm}

\noindent(1.2) $u$ a au moins un vecteur d'orbite dense

\vspace{2mm}

\noindent vérifient également la propriété (1.1) (voir le theorem 1.18 de \cite{BaMa}). Signalons au passage que les opérateurs qui vérifient (1.2) sont appelés hypercycliques et que le livre \cite{BaMa} de Bayart et Matheron constitue l'ouvrage de référence pour leur étude générale.

Le cas $r=1$ de notre théorème est déjà connu ; il découle du théorème 5.4 de \cite{CoMaPa}. Donnons quelques détails.
On choisit un opérateur diagonal de $\ell^2(\N^*)$ dont l'ensemble des éléments diagonaux est l'ensemble $e^{i\pi\Q}$ des racines de l'unité. On voit qu'un tel opérateur est rigide avec la suite des factorielles comme suite de rigidité. De plus on sait (problem 63 de \cite{Hal}) que son spectre est l'adhérence de l'ensemble de ses éléments diagonaux, en l'occurence $\overline{e^{i\pi\Q}}=\U$.

De manière générale, le spectre ponctuel d'un opérateur diagonal est égal à l'ensemble de ses éléments diagonaux. Par ailleurs il est immédiat que le spectre ponctuel d'un opérateur rigide est inclus dans $\U$. On ne peut donc pas puiser dans les opérateurs diagonaux pour construire les opérateurs
de notre théorème quand $0\le r<1$. Nos exemples montrent que cela devient possible si on élargit aux opérateurs diagonaux par blocs.

\section{Endomorphismes et opérateurs diagonaux par blocs}

Notons $\delta_k$ l'élément de $\C^{\N^*}$ défini par
$$
\delta_k(n) = \left|
\begin{array}{cll}
1 &\mathrm{si} &n=k\\
0 &\mathrm{si} &n\neq k.
\end{array}\right.
$$
On notera toujours $\ell^2:=\ell^2(\N^*)$ et $\ell^\infty:=\ell^\infty(\N^*)$. La famille $(\delta_k)_{k\ge1}$ est la base orthonormale canonique de $\ell^2$.

Rappelons que pour tout élément $a$ de $\ell^\infty$, la formule
$$
(D_ax)(k) = a(k)x(k)
$$
définit un opérateur $D_a$ de $\ell^2$ vérifiant $\|D_a\|_{\ell^2} = \|a\|_{\ell^\infty}$. On l'appelle l'opérateur diagonal de $\ell^2$ associé à $a$. On choisit le livre très pédagogique \cite{Hal} comme référence pour les opérateurs diagonaux. Plutôt qu'enchaîner définitions, lemmes et théorèmes, Halmos a dans \cite{Hal} la démarche de permettre au lecteur d'être le principal acteur de son initiation aux Hilbert et à leurs morphismes. L'objet du livre est de proposer 250 exercices classés par thème où alternent l'étude d'exemples bien choisis et des  éléments de théorie, en accompagnant tout cela par des commentaires éclairants.

Nous étendons pour notre part la famille des opérateurs diagonaux de $\ell^2$, à ce que nous appellerons la famille des opérateurs diagonaux par blocs de $\ell^2$. On définit ceux--ci comme des restrictions à $\ell^2$ d'endomorphismes de $\C^{\N^*}$.

Soit $(n(k))_{k\ge1}$ une suite d'entiers de $\N^*$. On note $E_k=\C^{n(k)}$. On \og découpe \fg{} $E=\C^{\N^*}$ en $\C^{n(1)+n(2)+\cdots+n(k)+\cdots} \simeq \prod\limits_{k=1}^{+\infty} \C^{n(k)}= \prod\limits_{k=1}^{+\infty}E_k$. On note $y=(y_k)_{k\ge1}\in \prod\limits_{k=1}^{+\infty} E_k$ un vecteur générique de $\C^{\N^*}$.

Pour alléger les notations, on écrira dorénavant \og $k$ \fg{} à la place de\break \og $k\ge1$ \fg{}.

On choisit pour tout $k$ un endomorphisme $u_k$ de $\C^{n(k)}$. On dira de l'endomorphisme $u=(u_k)_k$ de $\C^{\N^*}$ défini par
$$
u((y_k)_k) = (u_k(y_k))_k
$$
qu'il est un endomorphisme diagonal par blocs de $\C^{\N^*}$. 

\noindent On munit chaque espace $E_k=\C^{n(k)}$ de sa structure hermitienne usuelle. On utilise le symbole $\|\ \|$, à la fois pour les normes hilbertiennes des $E_k$, de $\ell^2$, et de leurs opérateurs. On a en particulier pour tout $y$ de $\ell^2$ 
$$
\|y\|^2 = \sum_k \|y_k\|^2.
$$
Dans le lemme suivant, le cas particulier des opérateurs diagonaux correspond au cas où tous les sous--espaces $E_k$ sont de dimension~1. L'énoncé et la preuve de ce lemme généralise les réponses aux questions posées aux problèmes 61 et 62 de \cite{Hal} relatifs aux opérateurs diagonaux.

\vspace{2mm}

{\monlem 

 Avec les endomorphismes $u_k$ et $u=(u_k)_k$ définis ci--dessus, les conditions suivantes sont équivalentes.

\begin{tabular}{ll}
\kern-0.7cm$(i)$ &$u(\ell^2) \subset \ell^2$\\

\kern-0.7cm$(ii)$ &$\sup\limits_k \|u_k\| <+\infty$.
\end{tabular}
}

\noindent\textit{Quand ces conditions sont réalisées, on a aussi}
$$
\kern-7cm \left| 
\begin{array}{ll}
&u_{|\ell^2} \ est\ continu\\
et\\
&\|u_{|\ell^2}\| = \sup\limits_k \|u_k\|
\end{array}
\right.\leqno(iii)
$$
\textit{et on dira que $u_{|\ell^2}$ est l'opérateur diagonal par blocs défini à partir des endomorphismes $u_k$.}

\vspace{2mm}

\textbf{Démonstration}

\vspace{2mm}

\noindent \fbox{$(i) \Rightarrow (ii)$} On raisonne par l'absurde ; on suppose donc que $\sup\limits_k \|u_k\|=+\infty$.  

\noindent On extrait une sous--suite $k_n$ pour laquelle $\|u_{k_n}\|\ge n$ pour tout $n$. En choisissant $y=(y_k)_k$ avec $\|y_k\|=1/n$ et $\|u_k(y_k)\| =\|u_k\|\cdot \|y_k\|$ quand $k=k_n$, et $y_k=0$ sinon, on voit que $y\in \ell^2$ et $u(y)\notin\ell^2$. Donc $u(\ell^2)\not\subset\ell^2$.

\vspace{2mm}

\noindent\fbox{$(ii) \Rightarrow [(i)$ et $(iii)]$}  Soit $y\in \ell^2$. On a 
$$
\sum_k \|u(y_k)\|^2 \le \sum_k \|u_k\|^2 \cdot \|y_k\|^2 \le \sup_k \|u_k\|^2 \cdot \|y\|^2.
$$
Donc avec (ii), $u(\ell^2) \subset \ell^2$, $u_{|\ell^2}$ est continu et
$$
\|u_{|\ell^2}\| \le \sup_k \|u_k\|.
$$ 
Par ailleurs, en considérant, pour chaque $k$, un vecteur de $\ell^2$ dont toutes les composantes sont nulles sauf la $k$-ième $y_k$ qui vérifie
 $\|y_k\| =1$ et $\|u_k(y_k)\|=\|u_k\|$, on a
$$
\|u_{|\ell^2}\| \ge \sup_k \|u_k (y_k)\| = \sup_k \|u_k\|.
$$

\vspace{2mm}

{\monlem Soit $u=(u_k)_k$ un endomorphisme diagonal par blocs de $\C^{\N^*}$ tel que
$$
\sup_k \|u_k\| <+\infty.
$$
Les conditions suivantes sont équivalentes

\begin{tabular}{ll}
\kern-0.7cm$(i)$ &$u_{|\ell^2}$ \textit{est inversible dans $\ell^2$}\\

\kern-0.7cm$(ii)$ &\textit{pour tout} $k$, $u_k$ \textit{est inversible, et}\\
\end{tabular}

$$
\sup_k \|u_k^{-1}\| <+\infty.
$$
 Quand ces conditions sont réalisées, on a aussi}
\vspace{1mm}

\begin{tabular}{ll}
\kern-0.7cm$(iii)$ &$\|(u_{|\ell^2})^{-1}\| = \sup_k \|u_k^{-1}\|$.
\end{tabular}

\vspace{2mm}

\textbf{Démonstration.}

\vspace{2mm}

\fbox{$(i) \Rightarrow [(ii)$ et $(iii)]$} Supposons $(i)$. On a alors pour tout $k$
$$
\ker u_k \subset (\prod_{j\ge 1} \ker u_j) \cap \ell^2 = (\ker u) \cap \ell^2 = \ker u_{|\ell^2}=\{0\}.
$$
Donc $u_k$ est inversible.

On en déduit que $u$ est inversible et
\begin{equation}
u^{-1} = (u_k^{-1})_k.
\end{equation}

Soit maintenant $y\in\ell^2$. Notons $x=(u_{|\ell^2})^{-1}(y)$. On a $u^{-1}(y) =u^{-1}(u\ld (x)) =x\in\ell^2$. Donc $u^{-1}(\ell^2)\subset \ell^2$ et $u^{-1}{\ld}=(u{\ld)}^{-1}$. En appliquant le lemme 2.1 à $u^{-1}$, on a donc avec (2.1) que $\sup\limits_k \|u_k^{-1}\| <+\infty$, $u^{-1} {}\ld$ est continu et
$$
\|(u\ld)^{-1}\| = \|u^{-1}{}\ld\| = \sup_k \|u_k^{-1}\|.
$$

\vspace{2mm}

\fbox{$(ii)\Rightarrow(i)$} Supposons $(ii)$. On a alors $u$ inversible et (2.1) est vérifié. En appliquant le lemme 2.1 à $u$, on a $u(\ell^2) \subset \ell^2$ et $u\ld$ est continu. En l'appliquant  cette fois à $u^{-1}$, on a $u^{-1}(\ell^2) \subset \ell^2$ et $u^{-1}{}\ld$ est continu. Il découle de tout cela que $u\ld$ est inversible d'inverse $u^{-1}{}\ld$.

On déduit immédiatement du lemme 2.2 le critère suivant d'appartenance au spectre.

\vspace{2mm}

{\monlem  Soit $u=(u_k)_k$ un endomorphisme diagonal par blocs de $\C^{\N^*}$ tel que
$$
\sup_k \|u_k\| <+\infty.
$$
Soit $\lambda\in \C$.

Les conditions suivantes sont équivalentes.

\begin{tabular}{ll}
\kern-0.7cm$(i)$& $\lambda\in \sigma (u\ld)$\\
\end{tabular}

$$
\kern-2.5cm \left|
\begin{array}{ll}
&\lambda
\in \bigcup\limits_k \sigma(u_k)\\
ou\\
&\lambda\notin \bigcup\limits_k \sigma(u_k)\ et\ \sup\limits_k \|(u_k-\lambda id)^{-1}\| =+\infty.
\end{array}.
\right.\leqno(ii)
$$}

\vspace{2mm}

\textbf{Remarque 1.} Dans le cas particulier des opérateurs diagonaux, le critère $(ii)$ s'écrit de manière plus simple
$$
\lambda\in \overline{\bigcup_k \sigma(u_k)}
$$
(voir le problème 63 de \cite{Hal}).

\vspace{2mm}

\textbf{Remarque 2.} Lorsque $u$ n'est pas pas inversible, écrivons par convention $\|u^{-1}\|=+\infty$. On peut alors écrire $(ii)$ sous la forme plus concise
$$
 \sup_k \|(u_k -\lambda id)^{-1}\| = +\infty.\leqno(ii)'
$$

\section{Notations}

Soient $f$ et $g$ deux fonctions définies sur une partie $D$ de $\R^d$ à valeurs dans $\R^{+*}$. On note $f\ll g$ pour signifier qu'il existe un réel $K\ge 1$ tel que pour tout élément $x$ de $D$, $f(x)\le Kg(x)$. On écrit $f\newsmile g$ quand on a simultanément  $f\ll g$ et $g\ll f$. Si $K$ dépend d'un paramètre $\eta$, c'est--à--dire est une fonction de $\eta$, on écrira par exemple $f\ll_\eta g$ à la place de $f\ll g$.

Soient $n$ et $p$ des entiers et $\alpha$ un réel tels que
\setcounter{equation}{0}
\begin{equation}
 1 \le p <n\qquad\mathrm{et}\qquad \alpha>1.
 \end{equation}
 Une fois choisis les entiers $n$ et $p$, on note $q$ l'entier défini par
 \begin{equation}
  p+q=n.
  \end{equation} 
On note
$$
\beta(n,p,\alpha)=(\beta_1(n,p,\alpha),\beta_2(n,p,\alpha),\ldots, \beta_n(n,p,\alpha)) = (\beta_1,\beta_2,\ldots,\beta_n)
$$
avec
\begin{equation}
\beta_j = \left|
\begin{array}{cl}
\alpha
 &(1\le j\le p)\\
 \noalign{\vskip2mm}
\dfrac{1}{\alpha^{p/q}} &(p<j\le n), 
 \end{array}
\right.
\end{equation}
et on prolonge $\beta$ par $n$--périodicité en posant $\beta_{j+n}:=\beta_j$. On note $e_1,e_2,\ldots,e_n$ la base canonique de $\C^n$ que l'on prolonge aussi par $n$--périodicité en posant $e_{j+n}:=e_j$. On note
$$
S(\beta(n,p,\alpha),t) = t^{n-1} \sum_{j=0}^{n-1} \frac{\beta_1\beta_2\cdots\beta_j}{t^j}.
$$
On note $v=v_{n,p,\alpha}$ l'endomorphisme de $\C^n$ défini par
\begin{equation}
v(e_j)=\beta_j e_{j+1}, \qquad (1\le j\le n)
\end{equation}
et pour $\lambda\in \C^*$
$$
\psi_\lambda(v)=\lambda^{n-1} \sum_{j=0}^{n-1}\Big(\frac{v}{\lambda}\Big)^j.
$$

\section{Spectre}

{\monlem Soit $(n,p,\alpha)$ vérifiant $(3.1)$.

\noindent On a pour tous $j$ et $k$ vérifiant $1\le j\le n-1$ et $1\le k\le n$,}
$$
\beta_k\ \beta_{k+1}\cdots \beta_{k+j-1} \le \beta_1\beta_2\cdots \beta_j.
$$

\vspace{2mm}
\setcounter{equation}{0}
\textbf{Démonstration.} On a
\begin{equation}
\beta_1\ge \beta_2\ge \cdots\ge \beta_n.
\end{equation}
L'inégalité demandée est donc immédiate si $k+j-1\le n$. Si $k+j-1>n$, on a d'après (4.1)
$$
\beta_\ell \le \beta_{\ell-(n-j)}\qquad \mathrm{pour}\ k\le \ell\le n.
$$
D'où
$$
\begin{array}{cl}
&\beta_k\beta_{k+1}\cdots \beta_{k+j-1} =(\beta_k\beta_{k+1}\cdots \beta_n)(\beta_1\beta_2\cdots\beta_{k+j-1-n})\\
\noalign{\vskip2mm}
\le &(\beta_{k-n+j}\beta_{k-n+j+1} \cdots\beta_j)(\beta_1\beta_2\cdots\beta_{k+j-1-n})\\
\noalign{\vskip2mm}
=&\beta_1\beta_2\cdots\beta_j.
\end{array}
$$

Dorénavant, à chaque fois qu'un nombre complexe $\lambda$ sera introduit, on notera $\rho=|\lambda|$.

\vspace{2mm}
 
{\monlem Soient $(n,p,\alpha)$ vérifiant $(3.1)$, $\lambda\in \C^*$ et $1\le j\le n-1$. On~a}
\begin{eqnarray}
 &&le\ \textit{polynôme}\ minimal\ de\ v\ est\ X^n-1,\\
 \noalign{\vskip2mm}
& &(v-\lambda id)\psi_\lambda(v) = (1-\lambda^n)id,\\
\noalign{\vskip2mm}
 & &\|v^j\| = \|v^j(e_1)\| = \beta_1\beta_2\cdots \beta_j,\\
 \noalign{\vskip2mm}
&& \|\psi_\lambda(v)\| \le S(\beta(n,p,\alpha),\rho),\\
\noalign{\vskip2mm}
&&\|\psi_\lambda(v)(e_1)\|^2 = S(\beta^2(n,p,\alpha),\rho^2).
 \end{eqnarray}
 
\textbf{Démonstration.}

D'après (3.3) et (3.2), on a $\beta_1\beta_2\cdots \beta_n=1$. Avec (3.4), on en déduit\break que $v$ permute circulairement les vecteurs de la base $(e_1,\beta_1e_2,\break \beta_1\beta_2e_3,\ldots, \beta_1\beta_2\cdots \beta_{n-1}e_n)$ de $\C^n$. D'où (4.2).

On a
$$
v\psi_\lambda(v)=\lambda^n \sum_{j=1}^{n}(v/\lambda)^j = \lambda\psi_\lambda(v)+v^n -\lambda^n id = \lambda\psi_\lambda(v)+(1-\lambda^n)id
$$
d'après (4.2). Cela entraîne (4.3).

Avec (3.4) on a d'une part $\|v^j e_1\| =\beta_1\beta_2\cdots \beta_j$ d'où $\|v^j\| \ge \beta_1\beta_2\cdots\beta_j$, et d'autre part
\begin{eqnarray*}
\Big\|v^j\Big(\sum_{k=1}^n x_ke_k\Big)\Big\|^2 =\Big \| \sum_{k=1}^n x_k \beta_k\beta_{k+1}\cdots \beta_{k+j-1} e_{k+j}\Big\|^2\\
= \sum_{k=1}^n (\beta_k\beta_{k+1}\cdots \beta_{k+j-1}|x_k|)^2 \le (\beta_1\beta_2\cdots\beta_j)^2 \sum_{k=1}^n |x_k|^2,
\end{eqnarray*}
d'après le lemme 4.1. On a donc aussi $\|v^j\|\le \beta_1\beta_2\cdots \beta_j$. Cela achève la preuve de (4.4).

En utilisant (4.4) on a
$$
\begin{array}{c}
\|\psi_\lambda(v)\| \le \rho^{n-1} \displaystyle\sum_{j=0}^{n-1} \frac{\|v\|^j}{\rho^j}\\
\noalign{\vskip2mm}
=\rho^{n-1} \displaystyle\sum_{j=0}^{n-1} \frac{\beta_1\beta_2\cdots \beta_j}{\rho^j} = S(\beta(n,p,\alpha),\rho)
\end{array}
$$
ce qui montre (4.5).

Enfin la formule (4.6) découle de (3.4) et (4.4).

{\monlem Soient $\lambda\in \C$ et $A\in \R$. On suppose que $0<\rho<1<A$. On a alors pour $(n,p,\alpha)$ vérifiant $(3.1)$ et $\alpha\le A$,}
\begin{equation}
S(\beta(n,p,\alpha),\rho) \newsmile{\!}_{\rho,A}\  \alpha^p \rho^q \sum_{j=0}^{q-1}\Big(\frac{1}{\rho \alpha^{p/q}}\Big)^j
\end{equation}
\textit{et}
\begin{equation}
S(\beta^2(n,p,\alpha),\rho^2) \newsmile{\!}_{\rho,A}\  \alpha^{2p} \rho^{2q} \sum_{j=0}^{q-1} \Big(\frac{1}{\rho^2\alpha^{2p/q}}\Big)^j.
\end{equation}

\textbf{Démonstration.}

On a d'après (3.3)
$$
\sum_{j=0}^{p-1} \frac{\beta_1\beta_2\cdots\beta_j}{\rho^j} = \sum_{j=0}^{p-1} \Big(\frac{\alpha}{\rho}\Big)^j = \frac{(\alpha/\rho)^p-1}{\alpha/\rho-1} \newsmile{\!}_{\rho,A} (\alpha/\rho)^p
$$
et
$$
\sum_{j=p}^{n-1} \frac{\beta_1\beta_2\cdots\beta_j}{\rho^j} = (\alpha/\rho)^p \sum_{j=0}^{q-1} \Big(\frac{1}{\rho \alpha^{p/q}}\Big)^j.
$$
On en déduit (4.7). On montre de même (4.8).

\vspace{2mm}

{\monlem Soit $(n(k),p(k),\alpha(k))_{k\ge2}$ une famille d'éléments de $\N^{*2}\times ]1,+\infty[$ telle que pour tout $k\ge 2$, $p(k)<n(k)$, }
\begin{equation}
la\ suite\ (n(k))_{k\ge 1}\ n'est\ pas\ \textit{bornée}
\end{equation}
\textit{et}
\begin{equation}
\lim_{k\rg +\infty} \alpha(k)=1.
\end{equation}
\textit{Alors l'opérateur diagonal par blocs $u\ld$ de $\ell^2$ avec $u=(u_k)_{k\ge1}$ où $u_1=id_{|\C}$ et $u_k=v_{n(k),p(k),\alpha(k)}$ pour $k\ge2$, a pour spectre}
$$
\sigma(u\ld) = \{\lambda\in\C : r \le |\lambda|\le 1\}
$$
\textit{où}
$$
r=\sup\{ s\in [0,1] : \textit{\ la\ suite\ } (\alpha(k)^{p(k)} s^{n(k)-p(k)})_{k\ge 2}\ \textit{est\ bornée}\ \}.
$$

\textbf{Démonstration.}

Commençons par remarquer que d'après (4.2), on a
\begin{equation}
\bigcup_{k\ge1} \sigma(u_k) \subset \U.
\end{equation}
Soient $\lambda\in \C^*$ et $k\ge 2$.

\vspace{2mm}

\textbf{1$\up{er}$ cas. $\rho>1$}

En utilisant successivement les formules (4.3) et (4.4) du lemme 4.2, la définition (3.3) et l'hypothèse (4.10), on~a
$$
\begin{array}{cl}
&\|(u_k-\lambda id)^{-1}\| = \dfrac{\|\psi_\lambda(u_k)\|}{|1-\lambda^{n(k)}|} \ll_\rho\displaystyle \sum_{j=0}^{n(k)-1} \frac{\|v_{n(k),p(k),\alpha(k)}^j\|}{\rho^j}
\\
=&\displaystyle\sum_{j=0}^{n(k)-1} \frac{\beta_1(n(k),p(k),\alpha(k))\beta_2(n(k),p(k),\alpha(k))\cdots \beta_j(n(k),p(k),\alpha(k))}{\rho^j}
\\
\le &\displaystyle\sum_{j=0}^{n(k)-1} \Big(\frac{\alpha(k)}{\rho}\Big)^j \ll_\rho \ 1.
\end{array}
$$
Avec (4.11) et le lemme 2.3, on en déduit que 
\begin{equation}
\textrm{le\ rayon\ spectral\ de \ } u\ld \textrm{\ est\ inférieur\ ou\ égal\ à\ un.}
\end{equation}

Pour les 2$\up{\textrm{ièmes}}$ et 3$\up{\textrm{ièmes}}$ cas, on note
$$
\tilde{r} = \sup \{s\in [0,1] : \ \mathrm{la\ suite\ } (\alpha(k)^{p(k)} s^{n(k)-p(k)})_{k\ge2}\ \textrm{est\ bornée}\}.
$$

\textbf{2$\up{\textrm{ième}}$ cas.} $0<\rho<1$ et la suite $(\alpha(k)^{p(k)}\rho^{n(k)-p(k)})_{k\ge 2}$ n'est pas bornée.

En utilisant successivement les formules (4.3) et (4.6) du lemme 4.2 et l'estimation (4.8) du lemme 4.3, on a alors
$$
\begin{array}{c}
\|(u_k-\lambda id)^{-1}\| \newsmile{\!}_{\rho} \|\psi_\lambda(u_k)\| \ge \|\psi_\lambda(u_k)(e_1)\|\\
\noalign{\vskip2mm}
= S^{1/2}(\beta^2(n(k),p(k),\alpha(k)),\rho^2) \gg_\rho \alpha(k)^{p(k)} \rho^{n(k)-p(k)}.
\end{array}
$$
Donc $\sup\limits_{k\ge 1}\|(u_k-\lambda id)^{-1}\| =+\infty$.

\noindent Avec le lemme 2.3, on en déduit que
\begin{equation}
\tilde{r} < |\lambda| <1 \Longrightarrow \lambda \in \sigma(u).
\end{equation}

\textbf{3$\up{\textrm{ième}}$ cas.} Soit $\gamma\in ]0,1[$ tel que la suite $(\alpha(k)^{p(k)}\gamma^{n(k)-p(k)})_{k\ge2}$ est bornée.

Soit alors $\lambda\in \C$ tel que $0<\rho<\gamma$.

Il existe $M\ge 1$ et $\delta>0$ tels que si
\begin{equation}
n(k) -p(k) \ge M,
\end{equation}
alors
\begin{equation}
\alpha
(k)^{p(k)/(n(k)-p(k))} \rho = (\alpha(k)^{p(k)} \gamma^{n(k)-p(k)})^{1/(n(k)-p(k))}\ \frac{\rho}{\gamma}<1-\delta.
\end{equation}
En utilisant successivement les formules (4.3) et (4.5) du lemme 4.2 et l'estimation (4.7) du lemme 4.3, on a
$$
\begin{array}{c}
\|(u_k-\lambda id)^{-1}\| \newsmile{\!}_\rho \|\psi_\lambda(u_k)\| \le S(\beta(n(k),p(k),\alpha(k)),\rho)\\
\noalign{\vskip2mm}
\newsmile{\!}_\rho \ Q := \alpha(k)^{p(k)} \rho^{n(k)-p(k)} \displaystyle\sum_{j=0}^{n(k)-p(k)-1} \Big(\frac{1}{\rho \alpha(k)^{p(k)/(n(k)-p(k))}}\Big)^j.
\end{array}
$$
On a de plus $Q\ll_\rho 1$. En effet pour les $k$ tels que l'on a (4.14), cela résulte de (4.15). Et pour les $k$ tels que $n(k)-p(k)<M$, la somme dans $Q$ est bornée, ainsi que la suite $(\alpha(k)^{p(k)}\rho^{n(k)-p(k)})_{k\ge 2}$. On a donc finalement $\|(u_k-\lambda id)^{-1}\| \ll_\rho 1$. Avec (4.11) et le lemme 2.3, on en déduit que
\begin{equation}
0 < |\lambda| <\tilde{r} \Longrightarrow \lambda \notin \sigma(u).
\end{equation}

Par ailleurs le spectre ponctuel de $u$ est
$$
\sigma_p(u)= \bigcup_{k\ge 1} \sigma(u_k) = \bigcup_{k\ge 2} \{z\in \C :z^{n(k)}=1\}
$$
d'après (4.2).

D'après l'hypothèse (4.9), on a donc
\begin{equation}
\U \subset \sigma(u).
\end{equation}

On conclut la preuve du lemme 4.4 en combinant (4.12), (4.13), (4.16), (4.17) et la compacité de $\sigma(u)$.

\section{Rigidité}

{\monlem Soit $(p{(k),\alpha(k)})_{k\ge2}$ une famille d'éléments de $\N^*\times]1,+\infty[$ avec $(\alpha(k))_{k\ge 2}$ décroissante et vérifiant $\alpha(k+1)=1+O(\frac{1}{k!})$, et pour tout $k\ge 2$, $(k-1)! \le p(k)<k!$. Alors l'opérateur diagonal par blocs $u_{|\ell^2(\N^*)}$ avec $u_1=id_{|\C}$ et pour tout $k\ge 2$, $u_k=v_{k!,p(k),\alpha(k)}$ est rigide de suite de rigidité $l!$. }

\vspace{2mm}

\textbf{Démonstration.}

Soit $y=(y_k)_{k\ge1}\in \ell^2$ avec $y_k\in E_k =\C^{k!}$. On a alors en utilisant les différentes hypothèses, les formules (4.2) et (4.4) du lemme 4.2, et la définition~(3.3),
$$
\begin{array}{l}
\|(u^{\ell!}-id)y\|^2 = \displaystyle\sum_{k\ge \ell+1} \|(u_k^{\ell!}-id)y_k\|^2\\
\noalign{\vskip2mm}
\le \displaystyle\sum_{k\ge \ell+1}
\|u_k^{\ell!}-id\|^2 \|y_k\|^2 \le 4 \displaystyle\sum_{k\ge \ell+1} \|u_k^{\ell!}\|^2\ \|y_k\|^2\\
\noalign{\vskip2mm}
=4 \displaystyle\sum_{k\ge \ell+1} (\alpha(k))^{2\ell!} \|y_k\|^2 \le 4 (\alpha(\ell+1))^{2\ell!} 
\displaystyle\sum_{k\ge \ell+1} \|y_k\|^2\\
\noalign{\vskip2mm}
\ll\displaystyle\sum_{k\ge \ell+1} \|y_k\|^2.
\end{array}
$$
On a donc $\lim\limits_{\ell\rg +\infty} u^{\ell!}y=y$.

\section{Preuve du théorème}

Le cas $r=1$ a été traité dans l'introduction. Quand $0\le r<1$, on choisit l'opérateur diagonal par blocs $u\ld$ de $\ell^2$ où 
\vspace{2mm}

\noindent$u=(u_k)_{k\ge1}$ avec $u_1=id_{|\C}$,

\vspace{2mm}

\noindent $u_k=v_{k!,p(k),\alpha(k)}$ pour tout $k\ge 2$, et
$$
(p(k),\alpha(k)) = \left|
\begin{array}{rll}
\Big(\max(k!-k,1),1+\dfrac{\log(1/r)}{(k-1)!}\Big) &\mathrm{si}\ 0<r<1\\
\noalign{\vskip2mm}
\Big(k!-1,1+\dfrac{1}{(k-1)!}\Big) &\mathrm{si}\ r=0.
\end{array}\right.
$$
On conclut la preuve du théorème en combinant le lemme 4.4 et le lemme~5.1.

\vspace{4mm}

\section{Annexes}

On utilisera pour ces deux annexes la notation
$$
\overline{D(z_0,r)} = \{z\in \C : |z -z_0| \le r\}.
$$

\noindent\textbf{Annexe 1. Preuve alternative de la Proposition CMP.}

\vspace{2mm}
Soit $u$ un opérateur d'un Banach et $n\ge 1$ tels que
\setcounter{equation}{0}
\begin{equation}
\|u^n -id\| \le 1/2.
\end{equation}
On a alors
$$
\sigma(u^n-id) \subset \overline{D(0,1/2)}
$$
d'où
$$
\sigma(u^n) \subset \overline{D(1,1/2)} \subset \{z \in \C : 1/2 \le |z| \le 3/2\}
$$
et
\begin{equation}
\sigma (u) \subset \{z \in \C : (1/2)^{1/n} \le |z| \le (3/2)^{1/n}\}.
\end{equation}
Maintenant si $u$ est uniformément rigide, il existe des entiers $n$ arbitrairement grands vérifiant (7.1) et donc aussi (7.2). On en déduit que $\sigma(u) \subset \U$.

\eject

\noindent\textbf{Annexe 2.} 

La preuve de la seconde partie de la Proposition 2.18 de \cite{CoMaPa}, de Costakis, Manoussos et Parissis, consiste à remarquer que le résultat suivant découle immédiatement du Corollary 1.2 du travail \cite{Mul} de Müller de 2001.

\vspace{2mm}

\noindent \textsc{Théorème} M
\textit{Le rayon spectral d'un opérateur rigide est inférieur ou égal à~un.}

\vspace{2mm}

\noindent\textbf{Preuve alternative du théorème M}

Soit $u$ un opérateur rigide d'un Banach de suite de rigidité $(n_k)_{k\ge1}$. Alors pour tout élément  $x$ de $E$, la suite de vecteurs $(u^{n_k}x)_{k\ge 1}$ est bornée. D'après le théorème de Banach--Steinhaus, la suite des normes des opérateurs $u^{n_k}$ est également bornée, disons par $M$. Alors pour tout entier $k$ supérieur ou égal à un, on a
$$
\sigma(u)\subset \overline{D(0,M^{1/n_k})},
$$
et on conclut en faisant tendre $k$ vers l'infini.

\vspace{2mm}

\noindent\textbf{Remerciements}

\vspace{2mm}

Nous remercions Etienne Matheron de nous avoir signalé l'article \cite{CoMaPa} de Costakis, Manoussos et Parissis, et Sophie Grivaux de nous avoir informé du lien entre son travail \cite{GrRo} avec Maria Roginskaya et le présent article.

\eject

\vskip4mm

\begin{tabular}{lll}

\kern-0.7cmPierre Mazet &\hspace{3.5cm}Eric Saias \\

\kern-0.7cm \textsf{piermazet@laposte.net} &\hspace{3.5cm}Sorbonne Université\\

&\hspace{3.5cm}LPSM\\

&\hspace{3.5cm}4, place Jussieu\\

&\hspace{3.5cm}75252 Paris Cedex 05 (France)\\

\vspace{2mm}

 &\hspace{3.5cm}\textsf{eric.saias@upmc.fr}
\end{tabular}

\end{document}